\documentclass{article}
\usepackage{amsmath, epsfig,amssymb, multicol,tikz}
\usepackage{color}
\newtheorem{lemma}{Lemma}
\newtheorem{theorem}{Theorem}

\newtheorem{proposition}{Proposition}

\newcommand{\M}{{\cal M}}

\renewcommand{\l}{{\cal L}}

\newcommand{\f}{{\bf f}}
\newcommand{\x}{{\bf x}}
\newcommand{\B}{{\cal B}}
\newcommand{\ignore}[1]{}

\title{Signless Laplacian spectral radius and  Hamiltonicity of graphs with large minimum degree}
\author{Yawen Li \thanks{Center for Applied Mathematics, Tianjin University,  Tianjin,    300072 China,
({\tt yawenli@tju.edu.cn}), }
\and Yao Liu  \thanks{Center for Applied Mathematics, Tianjin University,  Tianjin,    300072 China,
({\tt liuyao@tju.edu.cn}),}
\and Xing Peng  \thanks{Center for Applied Mathematics, Tianjin University,  Tianjin,    300072 China,
({\tt x2peng@tju.edu.cn}), Research is supported in part by National Natural Science Foundation of China (No. 11601380).}
}
\date{}
\begin{document}
\maketitle
\begin{abstract}
In this paper,  we establish  a tight sufficient  condition for the Hamiltonicity of graphs with large minimum degree in terms of the  signless Laplacian spectral radius and   characterize all extremal graphs.    Moreover, we prove a similar   result for balanced bipartite graphs.   Additionally, we construct  infinitely many  graphs   to show that results proved in this paper give new strength  for one to  determine the Hamiltonicity of graphs.
% It is surprising that theorems  in this paper contain  larger families of extremal  graphs than    previous results using the spectral radius of the adjacency matrix.
\end{abstract}

\section{Introduction}
A Hamiltonian cycle in a graph is a cycle that visits all vertices.  It is well known that the problem of determining the Hamiltonicity of graphs is remarkably difficult \cite{Karp}.  Thus it is meaningful to find sufficient conditions for graphs to be Hamiltonian.  A seminal  result due to Dirac \cite{Dirac} states that if $\delta(G) > \tfrac{n}{2}$, then $G$ is Hamiltonian, here $\delta(G)$ is the minimum degree of $G$ and $n$ is the number of vertices of $G$.

While many known results are in terms of  vertex degrees and the number of edges,  Krivelevich and Sudakov  showed the following  breakthrough  result using eigenvalues of the adjacency matrix.  Namely, they  proved that   a $d$-regular graph  is Hamiltonian if  its second largest eigenvalue (in the  absolute value) of  the adjacency matrix is sufficiently less than $d$.  Butler and Chung  \cite{BC} extended their result to general graphs by using eigenvalues of the combinatorial Laplacian matrix.  Both of   two results mentioned above used the Po\'sa's rotation method together with  some discrepancy and discrete
isoperimetric inequalities of graphs.

In 2009, Fiedler and Nikiforov \cite{fn} gave lower bounds on $\lambda(G)$ that imply the Hamiltonicity of $G$, where $\lambda(G)$ is the largest eigenvalue of the adjacency matrix of $G$.  After that,  many  other  spectral conditions  in the same spirit  for the Hamiltonicity of graphs  have been discovered, see e.g., \cite{be,ln1,llt,N16,ng,zhou}. In this paper, we will    first establish a  signless Laplacian  analogue of a  result due to Nikiforov \cite{N16} for graphs with large minimum degree. Then we strengthen   a result by Li and Ning \cite{ln1} for balanced bipartite graphs.  Before we state  our  theorems, we briefly recall results from \cite{N16}. We need to introduce two families of graphs.

 We use $K_n$ and $\overline{K}_n$ to denote the complete graph with $n$ vertices and the edgeless graph with $n$ vertices respectively. For two vertex disjoint graphs $G$ and $H$,  we write $G+H$ for their disjoint union which satisfies $V(G+H)=V(G)\cup V(H)$ and $E(G+H)=E(G)\cup  E(H)$.  We  write $G\vee H$ for their join  which satisfies $V(G\vee  H)=V(G)\cup V(H)$ and $E(G\vee H)=E(G) \cup E(H) \cup \{xy:x\in V(G), y\in V(H)\}$.

First,  for any $k >  1$ and $n\geq 2k+1$, let
 \[
M_k(n)=K_k\vee (K_{n-2k}+\overline{K}_k).
\]
The graph $M_k(n)$ is obtained from $K_{n-k}$ and ${\overline K}_k $ by connecting all vertices of ${\overline K}_k $  to all vertices of a $k$-subset of $K_{n-k}$.

Second, for any $k\geq1$ and $n\geq k+2$, let
\[
L_k(n)=K_1\vee(K_{n-k-1}+K_k).
\]
The graph $L_k(n)$ is obtained from $K_{n-k}$ and $K_{k+1}$ by identifying a vertex.   We note for any admissible $k$ and $n$, graphs     $L_k(n)$  and  $M_k(n)$ are not Hamiltonian.
 Strengthening a result by Li and Ning \cite{ln1}, Nikiforov \cite{N16} proved the following theorem.
\begin{theorem} \label{Ni16}
Let $k >1$, $n \geq k^3+k+4$, and let $G$ be a graph of order $n$, with minimum degree $\delta(G) \geq k$. If
\[
\lambda(G) \geq n-k-1,
\]
then $G$ has a Hamiltonian cycle unless $G=M_k(n)$ or $G=L_k(n)$.
\end{theorem}
For a graph $G$, let $Q(G)=D+A$ be the signless Laplacian matrix, here $D$ is the diagonal matrix of degrees.  We use $q(G)$ to denote the largest eigenvalue of $Q(G)$.
 If  $E' \subseteq E(G)$, then  we will use $G-E'$ to denote the subgraph of $G$ by deleting edges from $E'$. We shall prove a signless Laplacian analogue of Theorem \ref{Ni16}. To state our theorem, we need to define a collection of subgraphs  of  $M_k(n)$  and  $L_k(n)$.

For the graph $M_k(n)$,   let $X=\{v \in V(M_k(n)) : d(v)=k\}$,  $Y=\{v \in V(M_k(n)) : d(v)=n-1\}$, and $Z=\{v \in V(M_k(n)) : d(v)=n-k-1\}$.  Let  $E_1(M_k(n))$ be the set of  those edges of $M_k(n)$ whose both endpoints are from $Y \cup Z$.   We define
\[
\M_1(n,k)=\left\{G \subseteq M_k(n):  G=M_k(n)-E', \textrm{ where }  E' \subset E_1(M_k(n)) \textrm{ with }  |E'| \leq \left \lfloor \frac{k^2}{4} \right\rfloor \right\}.
\]

 Similarly, for the graph $L_k(n)$, we let $X=\{v \in V(L_k(n)) : d(v)=k\}$,  $Y=\{v \in V(L_k(n)) : d(v)=n-1\}$, and $Z=\{v \in V(L_k(n)) : d(v)=n-k-1\}$.  It is clear that the set $Y$ contains only one point.  We use  $E_1(L_k(n))$ to denote  edges of $L_k(n)$ whose both endpoints are from $Y \cup Z$. We define
\[
\l_1(n,k)=\left\{G \subseteq  L_k(n):  G=L_k(n)-E', \textrm{ where }  E' \subset E_1(L_k(n)) \textrm{ with }  |E'| \leq \left \lfloor \frac{k}{4} \right\rfloor \right\}.
\]
Our first result is the following theorem.
\begin{theorem} \label{t:thm1}
Assume $k>1$ and $n\geq k^4+k^3+4k^2+k+6$.  Let $G$ be a connected graph with  $n$ vertices and  minimum degree $\delta(G)\geq k$. If
\[
q(G) \geq 2(n-k-1),
\]
then $G$ has a Hamilton cycle unless $G \in \M_1(n,k)$ or  $G \in \l_1(n,k)$.
\end{theorem}
%We note that Theorem \ref{t:thm1} has a larger set of extremal graphs than Theorem \ref{Ni16} does.  This is uncommon in this area of research.  To compare Theorem \ref{Ni16} and Theorem \ref{t:thm1},
To show Theorem \ref{t:thm1}  does not covered by Theorem \ref{Ni16},
 in the last section, we will construct infinitely many Hamiltonian graphs that satisfy the condition in Theorem \ref{t:thm1} but do not satisfy the condition in Theorem \ref{Ni16}.
The proof of Theorem \ref{t:thm1} relies on techniques in \cite{ln1} and the following two lemmas.
\begin{lemma} \label{l:m1}
Assume $k > 1$ and $n\geq k^4+k^3+4k^2+k+6$. For each $G \in \M_1(n,k) \cup \l_1(n,k)$,  we have $q(G)  \geq  2(n-k-1)$.
\end{lemma}
We will need another family of subgraphs of $M_k(n)$  and $L_k(n)$  defined as follows:
\[
\M_2(n,k)=\left\{G \subset M_k(n):  G=M_k(n)-E', \textrm{ where }  E' \subset E_1(M_k(n)) \textrm{ with }  |E'|=\left\lfloor \frac{k^2}{4} \right\rfloor +1\right\}.
\]
\[
\l_2(n,k)=\left\{G \subset L_k(n):  G=L_k(n)-E', \textrm{ where }  E' \subset E_1(L_k(n)) \textrm{ with }  |E'|=\left\lfloor \frac{k}{4} \right\rfloor +1\right\}.
\]
We have the following lemma.
\begin{lemma}\label{l:m2}
Assume $k > 1$ and $n\geq k^4+k^3+4k^2+k+6$.  For each $G \in \M_2(n,k) \cup \l_2(n,k)$,  we have $q(G)<2(n-k-1)$.
\end{lemma}
In order to state our second theorem on balanced bipartite graphs, we need to introduce one more family of graphs.  A bipartite graph is called balanced if its vertex parts  have the same size.  For $k >1$ and  $n \geq 2k$,   we write $B_k(n)$ for the graph  obtained from $K_{n,n}$ by deleting all edges in a subgraph $K_{k,n-k}$. We note $B_k(n)$ is not Hamiltonian. Li and Ning \cite{ln1} proved the following theorem.
\begin{theorem} \label{lin1}
Let $G$ be a balanced bipartite graph of order $2n$ and of minimum degree $\delta(G) \geq k >1$.
\begin{enumerate}
\item If $n \geq (k+1)^2$ and $\lambda(G) \geq \lambda(B_k(n))$, then $G$ is Hamiltonian unless $G=B_k(n)$.
\item If $n \geq (k+1)^2$ and $q(G) \geq  q(B_k(n))$, then $G$ is Hamiltonian unless $G=B_k(n)$.
\end{enumerate}
\end{theorem}
For the graph $B_k(n)$,   let $S$ and $T$ be the vertex parts such that  the degree of vertices from $T$ is either $n$ or $n-k$.  Let $X=\{v \in  S : d(v)=k\}$,  $Y=\{v \in   T : d(v)=n\}$,   $W=\{v \in   T : d(v)=n-k\}$,
and $Z=\{v \in  S : d(v)=n\}$.  We note $S=X \cup Z$ and $T=Y \cup W$. We define $E_1(B_k(n))$ as those edges of $B_k(n)$ whose both endpoints are from $Y \cup W \cup Z$.   Let
\[
\B_1(n,k)=\left\{G \subseteq B_k(n):  G=B_k(n)-E', \textrm{ where }  E' \subseteq E_1(B_k(n)) \textrm{ with }  |E'| \leq \left \lfloor \frac{k^2}{4} \right\rfloor \right\}.
\]
We strengthen Part 2 of Theorem \ref{lin1} as follows.
\begin{theorem} \label{t:thm2}
Assume $k> 1$ and $n\geq k^4+3k^3+5k^2+5k+4$.  Let $G$ be a balanced bipartite graph with  $2n$ vertices and  minimum degree $\delta(G)\geq k$. If
\[
q(G) \geq 2n-k,
\]
then $G$ has a Hamilton cycle unless $G \in \B_1(n,k)$.
\end{theorem}
Here we notice $q(K_{n,n-k})=2n-k$ and $q(B_k(n))>2n-k$ since it contains $K_{n,n-k}$ as a proper subgraph. Thus the condition in Theorem \ref{t:thm2} is weaker than the one in Part 2 of Theorem \ref{lin1}. This is the reason why Theorem \ref{t:thm2} involves a  larger family of exception graphs  than Part 2 of Theorem \ref{lin1} does.   To show  Theorem \ref{t:thm2} is not covered by Part 1 of Theorem \ref{lin1}, in the last section,   we will construct infinitely many Hamiltonian graphs that satisfy the condition in Theorem \ref{t:thm2} but do not satisfy the condition in Part 1 of Theorem \ref{lin1}.
Similar to the proof of Theorem \ref{t:thm1}, we need another family of subgraphs of $B_k(n)$ defined as follows.
\[
\B_2(n,k)=\left\{G \subset B_k(n):  G=B_k(n)-E', \textrm{ where }  E' \subset E_1(B_k(n)) \textrm{ with }  |E'| = \left \lfloor \frac{k^2}{4} \right\rfloor+1 \right\}.
\]
Ideas from \cite{ln1} together with the following  lemma prove Theorem \ref{t:thm2}.
\begin{lemma}\label{l:m3}
 Assume $k > 1$ and $n\geq k^4+3k^3+5k^2+5k+4$.
 \begin{enumerate}
\item   For each $G \in \B_1(n,k)$,  we have $q(G) \geq 2n-k$.
\item  For each $G \in \B_2(n,k)$, we have $q(G)  < 2n-k$.
\end{enumerate}
\end{lemma}
The paper is organized as follows. In section 2, we will introduce some notation and present necessary preliminary results. We will prove Theorem \ref{t:thm1} in section 3 and sketch the proof of Theorem \ref{t:thm2} in section 4.   In section 5, we will  construct graphs to show theorems proved in this paper are new and give a few concluding remarks.

\section{Notation and preliminaries}
All graphs in this paper are simple and finite.  For those notation not defined here, we refer readers to the monograph written by West \cite{west}. For a graph $G$ and $v \in V(G)$, let $d_G(v)$ be the degree of   $v$ and $N_G(v)$ be the neighborhood of $v$, i.e.,  $N_G(v)=\{u \in V(G): \{u,v\} \in E(G)\}$ and $d_G(v)=|N_G(v)|$.  If the graph $G$ is clear under the context, we will drop the subscript $G$.  For a subset $X \subseteq V(G)$,  let $G[X]$ be the subgraph of $G$ induced by $X$. We use $A(G)$ to denote the adjacency matrix of $G$.  The signless Laplacian matrix $Q(G)$ associated with $G$ is defined as $D+A$, here $D$ is the diagonal matrix of degrees.
If $\x$ is a column vector of size $|V(G)|$, then
\[
\langle  Q(G) {\x},  \x\rangle=\sum_{v \in V(G)} d(v) \x_v^2 + 2\sum_{\{u,v\} \in E(G)} \x_u \x_v.
\]
If $q(G)$ is the largest eigenvalue of $Q(G)$,  then  by  Rayleigh's principle we have
\[
q(G)=\max_{\x} \frac{\langle  Q(G) {\x},  \x\rangle}{\langle  {\x},  \x\rangle}.
\]
 Let $\bf f$ be the eigenvector corresponding to $q(G)$, i.e.,
\[
Q(G) {\bf f}=q(G) {\bf f}.
\]
By the famous Perron-Frobenius theorem \cite{gr}, we get $\f_v>0$ for each $v \in V(G)$ provided $G$ is connected.
By taking the $v$-entry of both sides and rearranging terms, we get
\begin{equation} \label{s2:eq1}
(q(G)-d(v))  \f_v= \sum_{u \sim v} \f_u.
\end{equation}
It is easy to show the following proposition.
\begin{proposition} \label{s2:p1}
For any $u, v \in V(G)$,  we have
\[
(q(G)-d(u))(\f_u-\f_v)=(d(u)-d(v))\f_v+\sum_{s \in N(u)\setminus N(v)} \f_s -  \sum_{t \in N(v)\setminus N(u)} \f_t.
\]
\end{proposition}
\noindent
{\bf Proof:} By recalling \eqref{s2:eq1}, we get
\begin{align*}
(q(G)-d(u)) \f_u&= \sum_{a \in N(u)} \f_a, \\
(q(G)-d(v)) \f_v &= \sum_{b \in N(v)} \f_b.
\end{align*}
Therefore,
\begin{align*}
(q(G)-d(u))(\f_u-\f_v)&= (q(G)-d(u))\f_u-(q(G)-d(v))\f_v+(d(u)-d(v))\f_v \\
                                      &=(d(u)-d(v))\f_v+\sum_{s \in N(u)} \f_s -  \sum_{t \in N(v)} \f_t \\
                                      &=(d(u)-d(v))\f_v+\sum_{s \in N(u)\setminus N(v)} \f_s -  \sum_{t \in N(v)\setminus N(u)} \f_t. \\
\end{align*}
The proposition is proved.  \hfill $\square$

We will  repeatedly use the proposition above.  The following theorem from \cite{fy} provides an upper bound for $q(G)$.
\begin{theorem} \label{fy1}
Let $G$ be a graph of order $n$. Then
\[
q(G)\leq \frac{2e(G)}{n-1}+n-2.
\]
\end{theorem}
For balanced bipartite graphs, we will apply the following theorem from \cite{ln1}.
\begin{theorem}\label{t:ln1}
Let $G$ be a balanced bipartite graph of order $2n$. Then
\[
q(G)\leq \frac{e(G)}{n}+n.
\]
\end{theorem}
In the course of proving our results, we will need the following two theorems from \cite{ln1}.
\begin{theorem}\label{t:ln2}
Let $G$ be a graph of order $n\geq6k+5$, where $k\geq1$. If $\delta(G)\geq k$ and
\[
e(G)>\binom{n-k-1}{2}+(k+1)^2,
\]
then $G$ is Hamiltonian unless $G\subseteq L_k(n)$ or $G\subseteq M_k(n)$.
\end{theorem}
\begin{theorem}\label{t:ln3}
Let $G$ be a balanced bipartite graph of order $2n$. If $\delta(G)\geq k\geq1,n\geq 2k+1$ and
\[
e(G)>n(n-k-1)+(k+1)^2,
\]
then $G$ is Hamiltonian unless $G\subseteq B_k(n)$.
\end{theorem}
\section{Proof of Theorem \ref{t:thm1}}
We start with this section by proving Lemma \ref{l:m1}.

\noindent
{\bf Proof of Lemma \ref{l:m1}:}  Let  $G \in \M_1(n,k) \cup \l_1(n,k)$.  Recall three subsets $X, Y,$ and $Z$ of $V(G)$ for both cases. For each case,  we define a column vector $\bf h$ such that ${\bf h}_u=1$ for all $u \in Y\cup Z$ and ${\bf h}_v=0$ for all $v \in X$.  We  note  $q(\overline{K}_k+K_{n-k})=2(n-k-1)$ and $\bf h$ is the corresponding eigenvector.  If  $G \in \M_1(n,k) $, then we get
\[
\langle Q(G) {\bf h}, {\bf h} \rangle - \langle Q(\overline{K}_k+K_{n-k}) {\bf h}, {\bf h} \rangle=k^2-4|E'| \geq 0.
 \]
 Similarly, for   $G \in \l_1(n,k) $,   we have
\[
\langle Q(G) {\bf h}, {\bf h} \rangle - \langle Q(\overline{K}_k+K_{n-k}) {\bf h}, {\bf h} \rangle=k-4|E'| \geq 0.
 \]
By Rayleigh's principle,  we obtain $q(G) \geq 2(n-k-1)$ in each case.  \hfill $\square$

In order to prove Lemma \ref{l:m2}, we need a lower bound on $q(G)$.
\begin{proposition} \label{s3:p1}
For each $G \in \M_2(n,k) \cup \l_2(n,k)$, we have $q(G)>2n-2k-3$.
\end{proposition}
\noindent
{\bf Proof:} Suppose $G \in \M_2(n,k) \cup \l_2(n,k)$.  Let $\bf h$  be  the vector defined  in the proof of Lemma  \ref{l:m1}.   In the case of  $G \in \M_2(n,k)$, we have
\[
\langle Q(G) {\bf h}, {\bf h} \rangle - \langle Q(\overline{K}_k+K_{n-k}) {\bf h}, {\bf h} \rangle=k^2-4|E'| \geq -4.
\]
In the case of  $G \in \l_2(n,k)$, we have
\[
\langle Q(G) {\bf h}, {\bf h} \rangle - \langle Q(\overline{K}_k+K_{n-k}) {\bf h}, {\bf h} \rangle=k-4|E'| \geq -4.
\]
In each case, we have $q(G) \geq 2(n-k-1)-\tfrac{4}{\|\bf h\|^2} > 2n-2k-3$.  \hfill $\square$

We next prove Lemma \ref{l:m2}.  In the following propositions, we will assume $G \in \M_2(n,k)$ and give their detailed proofs.  Since arguments in  the case of $G \in \l_2(n,k)$ are similar to those in  the case of $G \in \M_2(n,k)$, we will only sketch them.

Let $G$ be a graph from $\M_2(n,k)$ with the maximum signless Laplacian spectral radius.  Moreover,  we assume  $G[Y]$ contains the largest number of edges. Let $\f$ be the eigenvector corresponding to $q(G)$. We assume further $\max_{v \in V(G)} \f_v =1$.   Our goal is to show the vector $\f$ is close (entrywise) to the vector $\bf h$ as defined in the proof of Lemma \ref{l:m1}.

We define two subsets of $Y $ as follows:
\[
Y_1=\{y \in Y: d(y)=n-1\}  \textrm{  and  } Y_2=\{y \in Y: d(y) \leq n-2\}.
\]
 Similarly, we define two subsets of $Z $ as follows:
\[
Z_1=\{z \in Z: d(z)=n-k-1\}  \textrm{  and  } Z_2=\{z \in Z:  d(z) \leq n-k-2\}.
\]
We note $Z_1 \not = \emptyset$ as $n \geq k^4+k^3+4k^2+k+6$.   To compare the difference between $\max_{v \in V(G)} \f_v $ and $\min_{v \in Y\cup Z} \f_v$, we need to prove the following  propositions.

\begin{proposition}\label{s3:p2}
Assume $G \in \M_2(n,k)$ as defined above. For each $x \in X$, we have $\f_x \leq \tfrac{k}{q(G)-k}$.
\end{proposition}
\noindent
{\bf Proof:} Applying equation \eqref{s2:eq1} with $x$, we get
\[
(q(G)-d(x)) \f_x= \sum_{y \in  Y} \f_y.
\]
Since $d(x)=k$ and  $\max_{v \in V(G)} \f_v=1$, the lemma follows. \hfill $\square$
\begin{proposition}\label{s3:p3}
Assume $G \in \M_2(n,k)$ as defined above.  If $Y_2 \not = \emptyset$, then we have  $\f_y  < \f_z$ for all $y \in Y_2$ and $z \in Z_1$.
\end{proposition}
\noindent
{\bf Proof:} We assume that there are some $y \in Y_2$ and $z \in Z_1$ such that $\f_y  \geq  \f_z$. Let $w$ be a vertex from $Y$ such that $\{y, w\}$ is a non-edge. We define a new graph  $G' \in \M_2(n,k)$ by  removing $\{w,z\}$ and adding $\{y, w\}$,  as shown in Figure 1.  Since
\[
\langle Q(G') {\f}, {\f} \rangle - \langle Q(G) {\f}, {\f} \rangle=(\f_y-\f_z)(\f_y+\f_z+2\f_w) \geq 0,
\]
we get $q(G')  \geq  q(G)$ and $G'[Y]$ has more edges than $G[Y]$, which is a contradiction to the choice of $G$.  The proposition is proved. \hfill $\square$
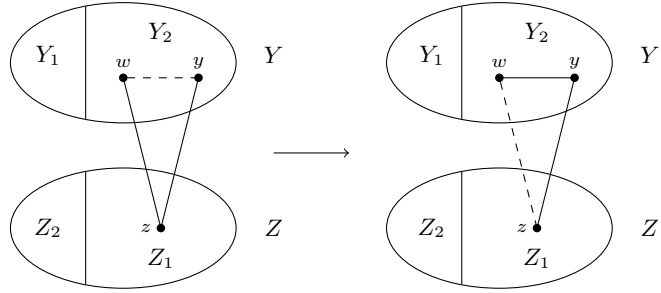
\begin{figure}[htp]
\begin{center}
\begin{tikzpicture}
\scriptsize
[inner sep=0pt]
%%%%%%%%%%%%%%%%%%%%%%%%%%%%%%%%%%%%%%%%%%%ÐÂÔö¶¥µã
\draw  (1,1.2) ellipse (1.5 and 0.8);
\draw  (-4,1.2) ellipse (1.5 and 0.8);
\draw  (1,-1) ellipse (1.5 and 0.8);
\draw  (-4,-1) ellipse (1.5 and 0.8);
\draw[->](-2,0)--(-1,0);
\filldraw [black] (1,1) circle (1.4pt)
                  (2,1) circle (1.4pt)
                  (1.5,-1) circle (1.4pt)
                   (-3,1) circle (1.4pt)
                  (-4,1) circle (1.4pt)
                  (-3.5,-1) circle (1.4pt);
\draw  (2,1)-- (1,1) ;\draw  (2,1)-- (1.5,-1) ;
\draw  (-3.5,-1)-- (-4,1) ;\draw  (-3,1)-- (-3.5,-1) ;
\draw[style=dashed]  (1.5,-1)-- (1,1) ;
\draw[style=dashed]   (-3,1)-- (-4,1) ;
\draw  (-3,1.2) node { $y$};
\draw   (-4,1.2)node { $w$};
\draw   (-3.7,-1) node { $z$};
\draw  (1,1.2) node { $w$};
\draw   (2,1.2)node { $y$};
\draw   (1.3,-1) node { $z$};
\draw   (-2,1.3) node { \small $Y$};
\draw  (-2,-1) node { \small $Z$};
\draw   (3,1.3) node { \small $Y$};
\draw  (3,-1) node { \small $Z$};
\draw  (0.5,1.95)-- (0.5,0.45) ;\draw  (0.5,-1.75)-- (0.5,-0.25) ;
\draw  (-4.5,1.95)-- (-4.5,0.45) ;\draw  (-4.5,-1.75)-- (-4.5,-0.25) ;
\draw   (-5,1.3) node { \small $Y_1$};
\draw  (-5,-1) node { \small $Z_2$};
\draw   (-3.5,1.6) node { \small $Y_2$};
\draw  (-3.5,-1.4) node { \small $Z_1$};
\draw   (0.1,1.3) node { \small $Y_1$};
\draw  (0.1,-1) node { \small $Z_2$};
\draw   (1.5,1.6) node { \small $Y_2$};
\draw  (1.5,-1.4) node { \small $Z_1$};
\end{tikzpicture}
\caption{An example for $G$ and $G'$}
\label{figure1}
\end{center}
\end{figure}

\begin{proposition}\label{s3:p4}
Assume $G \in \M_2(n,k)$ as defined above.  If $Z_2 \not = \emptyset$, then we have  $\f_w  < \f_z$ for any $w \in Z_2$ and $z \in Z_1$.
\end{proposition}
\noindent
{\bf Proof:} We first notice    $N(w) \setminus \{z\} \subset N(z) \setminus \{w\}$  for each $w \in Z_2$ and $z \in Z_1$.  Applying  Proposition \ref{s2:p1} with $u=z$ and $v=w$, we get
\[
(q(G)-d(z))(\f_z-\f_w)=(d(z)-d(w))\f_w+\f_w-\f_z+\sum_{s \in N(z)\setminus (N(w) \cup w)} \f_s.
\]
Equivalently,
\begin{equation} \label{lm2:eq2}
(q(G)-d(z)+1)(\f_z-\f_w)=(d(z)-d(w))\f_w+\sum_{s \in N(z)\setminus (N(w) \cup w)} \f_s.
\end{equation}
 We note $d(w)<d(z) $ and the lemma follows. \hfill $\square$

Similarly, we can use Proposition \ref{s2:p1}  to  show the following one.
\begin{proposition}\label{s3:p5}
Assume $G \in \M_2(n,k)$ as defined above.
\begin{enumerate}
\item If $Y_1, Y_2 \not = \emptyset$, then we have  ${\f}_s  < \f_t$ for any $s \in Y_2$ and $t \in Y_1$.
\item  If $Y_1 \not = \emptyset$, then we have $\f_z < \f_y$ for any $y \in Y_1$ and $z \in Z_1$.
\end{enumerate}
\end{proposition}

The key step for proving  Lemma \ref{l:m2} is to show the following proposition.
\begin{proposition} \label{s3:p6}
Assume $G \in \M_2(n,k)$ as defined above.  We have
\[
\max_{v \in V(G)} \f_v - \min_{u \in Y \cup Z} \f_u \leq  \frac{k^2+2k+6}{2(q(G)-n+1)}.
\]
\end{proposition}
\noindent
{\bf Proof:}
We have two cases depending on $Y_1$.
\begin{description}
\item[Case 1:] $Y_1  = \emptyset$.  By  Proposition \ref{s3:p3} and  Proposition \ref{s3:p4}, we get  $\max_{v \in V(G)} \f_v $ is attained by vertices from $Z_1$.  Take $z$ as an arbitrary vertex from $Z_1$
.  We observe that $z$ is adjacent to all other vertices in $Y \cup Z$.

 If  $w \in Z_2$, then  we have  $N(z)\setminus (N(w) \cup w)=\{k: k \in Y \cup Z \textrm{ and } k \not \sim w\}$. Thus
we get $d(z)-d(w) \leq  \lfloor \tfrac{k^2}{4} \rfloor+1$ and $|N(z)\setminus (N(w) \cup w)| \leq \lfloor \tfrac{k^2}{4} \rfloor+1$.
    Recalling \eqref{lm2:eq2}, we  get
    \begin{align*}
    (q(G)-d(z)+1)(\f_z-\f_w)&=(d(z)-d(w))\f_w+\sum_{s \in N(z)\setminus (N(w) \cup w)} \f_s\\
                                             & \leq \frac{k^2}{2}+2.
    \end{align*}
    Since $d(z)=n-k-1$, we get
    \[
    \f_z -\f_w \leq \frac{k^2+4}{2(q(G)-n+k+2)} <\frac{k^2+2k+6}{2(q(G)-n+1)}.
    \]
      If  $w \in  Y_2$,  then  we have  $N(z)\setminus (N(w) \cup w)=\{k: k \in Y \cup Z \textrm{ and } k \not \sim w\}$ and $N(w) \setminus (N(z) \cup z)=X$.   Thus      $|N(z)\setminus (N(w) \cup w)| \leq \lfloor \tfrac{k^2}{4} \rfloor+1$.  We also observe      $|d(z)-d(w)| \leq \lfloor \tfrac{k^2}{4} \rfloor+k+1$.  Applying Proposition \ref{s2:p1} with $u=z$ and $v=w$, we get
         \begin{align*}
    (q(G)-d(z))(\f_z-\f_w) &=(d(z)-d(w))\f_w+\f_w+\sum_{s \in N(z) \setminus (N(w)\cup w)} \f_s - \f_z- \sum_{t \in X} \f_t\\
                                         & \leq   \lfloor \frac{k^2}{4} \rfloor + k+1+\f_w+\lfloor\frac{k^2}{4} \rfloor+1-\f_z-\sum_{t \in X} \f_t.
   \end{align*}
      Equivalently,
     \[
      (q(G)-d(z)+1)(\f_z-\f_w) \leq    \frac{k^2}{2} + k+2.
      \]
     We get
    \[
    \f_z -\f_w \leq \frac{k^2+2k+4}{2(q(G)-n+k+2)}<\frac{k^2+2k+6}{2(q(G)-n+1)}.
    \]
\item[Case 2:]  $Y_1   \not = \emptyset$.  By  Proposition \ref{s3:p4} and Proposition \ref{s3:p5}, we get $\max_{v \in V(G)}  {\bf  f}_v $ is attained by vertices from $Y_1$.  Let $z$ be a vertex from $Y_1$. By repeating the argument in Case 1, we can show
\[
\f_z-\f_w \leq \frac{k^2+2k+6}{2(q(G)-n+2)}<\frac{k^2+2k+6}{2(q(G)-n+1)}
\]
for any $w \in Y_2 \cup Z$.
\end{description}
The proof is completed. \hfill $\square$

We next assume $G \in \l_2(n,k)$. Recalling three subsets $X$, $Y$, and $Z$ of $V(G)$, here we suppose $Y=\{y\}$.   We take such $G$ satisfying $q(G)$ is maximum and $d(y)$ is the largest.  Let $\f$ be the eigenvector corresponding to $q(G)$.  Suppose $\f$ satisfies $\max_{v \in V(G)} \f_v =1$.   We define
 \[
 Z_1=\{z \in Z : d(z)=n-k-1\} \textrm{  and  } Z_2=\{z \in Z : d(z) \leq n-k-2\}.
 \]
 We will need the following proposition for the case of $G \in \l_2(n,k)$.
%  We can show the following proposition by repeating arguments in Proposition \ref{s3:p2}, Proposition \ref{s3:p3}, Proposition \ref{s3:p4}, Proposition \ref{s3:p5}, and   Proposition \ref{s3:p6} .
\begin{proposition} \label{s3:p7}
Let $G$ be assumed as above.
\begin{enumerate}
\item For each $x  \in X$, we have $\f_x \leq \tfrac{k}{q(G)-k}$.
\item If $Z_2 \not = \emptyset$, then we have $\f_w< \f_z$ for each $w \in Z_2$ and each $z \in Z_1$.
\item If $d(y) \leq n-2$, then $\f_y  < \f_z$ for each $z \in Z_1$.
\item If $d(y) =n-1$, then $\f_y  > \f_z$ for each $z \in Z_1$.
\item  We have  $\max_{v \in V(G)} \f_v - \min_{u \in Y \cup Z} \f_u \leq  \frac{k^2+2k+6}{2(q(G)-n+1)}.$
\end{enumerate}
\end{proposition}

\noindent
{\bf Proof:} Since proofs of  parts of this proposition are very similar to those of  Propositions  \ref{s3:p2}, \ref{s3:p3}, \ref{s3:p4},  \ref{s3:p5}, \ref{s3:p6}, we only give the sketch here.  Parts 1, 2, 3, 4 use ideas from proofs of Propositions \ref{s3:p2}, \ref{s3:p4},  \ref{s3:p3}, \ref{s3:p6} respectively.   For Part 5, we have two cases depending on $d(y)$.

\begin{description}
\item[Case 1:] $d(y) \leq n-2$.  We note    $\max_{v \in V(G)} \f_v$ is achieved by vertices from $Z_1$.  Let $z$ be an arbitrary vertex from $Z_1$.

If $w \in Z_2$, then we have  $N(z)\setminus (N(w) \cup w)=\{k: k \in Y \cup Z \textrm{ and } k \not \sim w\}$. Therefore,  $|N(z) \setminus (N(w) \cup w)| \leq \lfloor \tfrac{k}{4} \rfloor+1$ and $d(z)-d(w) \leq \lfloor \tfrac{k}{4} \rfloor+1$.  By Proposition \ref{s2:p1}, we get
\[
    (q(G)-d(z))(\f_z-\f_w)=(d(z)-d(w))\f_w+\f_w+\sum_{s \in N(z)\setminus (N(w) \cup w)} \f_s -\f_z.
\]
Therefore,
\[
   (q(G)-d(z)+1)(\f_z-\f_w) \leq \frac{k}{2}+2.
\]
  If $w=y$, then we have    $N(z)\setminus (N(w) \cup w)=\{k: k \in Y \cup Z \textrm{ and } k \not \sim w\}$ and $N(w) \setminus (N(z) \cup z)=X$.
  We notice $|d(z)-d(w)| \leq \lfloor \tfrac{k}{4} \rfloor+k+1$ and  $|N(z) \setminus (N(w) \cup w)| \leq \lfloor \tfrac{k}{4} \rfloor+1$.
   By Proposition \ref{s2:p1}, we get
   \begin{align*}
    (q(G)-d(z))(\f_z-\f_w) &=(d(z)-d(w))\f_w+\f_w+\sum_{s \in N(z) \setminus (N(w)\cup w)} \f_s - \f_z- \sum_{t \in  X} \f_t\\
                                         & \leq   \lfloor \frac{k}{4} \rfloor + k+1+\f_w+\lfloor\frac{k}{4} \rfloor+1-\f_z-\sum_{t \in X} \f_t.
   \end{align*}
    Therefore,
    \[
       (q(G)-d(z)+1)(\f_z-\f_w) \leq   \frac{k}{2}+k+1.
    \]
In both subcases, we have  $\f_z-\f_w \leq \tfrac{k^2+2k+6}{2(q(G)-n+1)}$.
\item[Case 2:] $d(y)=n-1$. We note $\max_{v\in V(G)} \f_v=\f_y$.  We can use the argument above to show the desired upper bound on $\f_y-\f_z$ for all $z \in Z$.
\end{description}
The proposition is proved. \hfill $\square$

We are ready to prove Lemma \ref{l:m2}.

\noindent
{\bf Proof of Lemma \ref{l:m2}:} We first assume $G \in \M_2(n,k)$ such that $G$ has the largest signless Laplaican spectral radius and $G[Y]$ contains the largest number of edges.  Let $q(G)$ be the largest eigenvalue of $Q(G)$ and $\f$ be the eigenvector of $q(G)$.  Recalling Proposition \ref{s3:p2} and Proposition \ref{s3:p6},  we get
\begin{align*}
\langle Q(G)\f, \f \rangle - \langle Q(\overline{K}_k + K_{n-k})\f, \f \rangle &=\sum_{x \in X, y \in Y} (\f_x+\f_y)^2-\sum_{\{u,v\} \in E'} (\f_u+\f_v)^2\\
                                                                                                                                   & \leq k^2\left( 1+\frac{k}{q(G)-k}\right)^2-4|E'| \left(1-\frac{k^2+2k+6}{2(q(G)-n+1)}  \right)^2\\
                                                                                                                                   & < 0.
\end{align*}
Here we used $n \geq k^4+k^3+4k^2+k+6$ and $q(G) > 2n-2k-3$ by Proposition \ref{s3:p1}.  We note $q(\overline{K}_k + K_{n-k})=2(n-k-1) \geq {\langle Q(\overline{K}_k + K_{n-k})\f, \f \rangle}/{\langle \f, \f \rangle}$. Thus we proved $q(G)<2(n-k-1)$ in this case.

 When $G \in \l_2(n,k)$, we can assume $q(G)$ is the maximum and $d(y)$ is the largest.  Applying Proposition \ref{s3:p7} and repeating the argument for $G \in \M_2(n,k)$, we can prove $q(G)<2(n-k-1)$ easily.
\hfill $\square$

We are ready to prove Theorem \ref{t:thm1}.

\noindent
{\bf Proof of Theorem \ref{t:thm1}:}   By theorem \ref{fy1}, we get
\[
2(n-k-1) \leq q(G)\leq\frac{2e(G)}{n-1}+n-2.
\]
Since we assume  $n\geq k^4+k^3+4k^2+k+6$, we  get
\begin{eqnarray*}
 e(G) &>& \frac{(n-2k)(n-1)}{2} \\
  {} &=& \binom{n-k-1}{2}+\frac{2n-k^2-k-2}{2} \\
  {} &\geq& \binom{n-k-1}{2}+(k+1)^2.
\end{eqnarray*}

By Theorem \ref{t:ln2},  $G$ has a Hamilton cycle unless $G\subseteq L_k(n)$ or $ G \subseteq M_k(n)$. Together with Lemma \ref{l:m1} and Lemma \ref{l:m2}, we complete the proof. \hfill $\square$\\

\section{Proof of Theorem \ref{t:thm2}}

\noindent
We first observe the following:   if   Lemma \ref{l:m3} holds,  then the combination of  Theorem \ref{t:ln1} and Theorem \ref{t:ln3} yields  Theorem \ref{t:thm2}. We are left to prove     Lemma \ref{l:m3}.

We note $q(K_{n,n-k}+\overline{K_k})=2n-k$.  Let $\f$ be an eigenvector  corresponding to $2n-k$.   If we assume $\max_{v} \f_v=1$, then  we have $\f_u=1$ for  $d(u)=n$ and $\f_u=1-\tfrac{k}{n}$ for  $d(u)=n-k$.  Recalling the definition of  $\B_1(n,k)$ and repeating the proof for Lemma \ref{l:m1}, we can prove Part 1 of Lemma \ref{l:m3}.  Given Proposition \ref{s4:p2},  we can show Part 2 of Lemma \ref{l:m3} by the same argument as the one in the proof of  Lemma \ref{l:m2}.

 Let $G \in \B_2(n,k)$ such that  $q(G)$ is the maximum and $G[Y \cup Z]$ induces the largest number of edges. We assume the corresponding  eigenvector is $\f$ and  $\max_v \f_v=1$.
We define $Y_1=\{v \in Y: d(v)=n \}$,  $Y_2=\{v \in Y: d(v) \leq n-1 \}$,  $W_1=\{v \in W: d(v)=n -k\}$,  $W_2=\{v \in W: d(v) \leq n-k-1 \}$,   $Z_1=\{v \in Z: d(v)=n \}$, and  $Z_2=\{v \in Z: d(v) \leq n-1 \}$.

It remains to establish  the following proposition.
\begin{proposition} \label{s4:p2}
Let $G$ be assumed as above.
\begin{enumerate}
\item $2n-k-1 \leq q(G) \leq 2n-k+1$.
\item For each $x \in X$, we have $\f_x \leq \tfrac{k}{q(G)-k}$.
\item  If $Y_2 \not = \emptyset$, then we have  $\f_y  < \f_w$ for each $y \in Y_2$ and each $w \in W_1$.
\item  If $Y_1 \not = \emptyset$, then we have  $\f_w  < \f_y$ for each $y \in Y_1$ and each $w \in W$.
\item  If $Y_1,Y_2 \not = \emptyset$, then we have  $\f_t  < \f_y$ for each $t \in Y_2$ and each $y \in Y_1$.
\item If $W_2 \not =\emptyset$, then we have  $\f_s  < \f_t$ for each $s \in W_2$ and each $t \in W_1$.
\item  If $Z_2 \not = \emptyset$, then we have  $\f_t  < \f_z$ for each $t \in Z_2$ and each $z \in Z_1$.
\item  We have  $\max_{v \in V(G)}\ \f_v - \min_{u \in Y \cup Z \cup W}\ \f_u \leq  \tfrac{3k^2+8k+20}{4(q(G)-n)}.$
\end{enumerate}
\end{proposition}

\noindent
\textbf{Proof:} The upper bound in Part 1 follows from Theorem \ref{t:ln1} and the proof of the lower bound uses the same idea as the one in Proposition \ref{s3:p1}.  Equation \eqref{s2:eq1}  gives Part 2.
When we prove Part 3, we have to apply Proposition \ref{s2:p1} and the edge-switching idea as we did in the proof of Proposition \ref{s3:p3}. Part 4-Part 7 follow from Proposition \ref{s2:p1} straightforwardly.

 For the proof of part 8,  we consider the following cases.
 \begin{description}
 \item[Case 1:] $Y_1 = \emptyset$. We observe $\max_{v \in V(G)} \f_v$ is achieved by vertices from $W_1 \cup Z_1$.  Suppose  $\f_w=\max_{v \in V(G)} \f_v$ for some $w \in W_1$.  We note $|W_1|\geq n-\lfloor \tfrac{k^2}{4} \rfloor -1$ and $\f_w=\f_{w'}$ for all $w, w' \in W_1$.  Let $z \in Z$.  We have $d(z) \leq n$.  By equation \eqref{s2:eq1}, we get $(q(G)-d(z)) \f_z =\sum_{v: v \in N(z)} \f_v$. Obviously,  we have $(q(G)-d(z)) \f_z \geq \sum_{w \in W_1} \f_w$. Therefore,
     \[
     \f_z \geq \frac{4n-k^2-4}{4(q(G)-n)}.
      \]
     Let $y$ be an arbitrary vertex from $(Y \cup W) \setminus W_1$.  We can apply Proposition \ref{s2:p1} with $u=w$ and $v=y$ to show
     \[
      \f_z-\f_y \leq \frac{k^2+4k+4}{2(q(G)-n+1)}.
     \]
If   $\f_z=\max_{v \in V(G)} \f_v$ for some $z \in  Z_1$, then we can use  equation \eqref{s2:eq1} to show a lower bound on $\f_y$ for each $y \in Y \cup W$.  We  apply  Proposition \ref{s2:p1} to get a lower bound on $\f_{z'}$ for all $z' \in Z_2$.   Tedious calculation can confirm Part 8 in this case.
 \item[Case 2:]   $Y_1 \not = \emptyset$. We observe $\max_{v \in V(G)} \f_v$ is achieved by vertices from $Y_1 \cup Z_1$.  We can repeat the argument in Case 1 to show Part 8 in this case.
 \end{description}
We completed the proof of this proposition. \hfill $\square$

\section{Concluding remarks and examples}
In this paper, we proved a new sufficient spectral condition for the Hamiltonicity of graphs. The main work is to  prove a lower bond on the difference between the largest entry and the smallest entry (of a particular subset of vertices) of the principal eigenvector of a family of graphs.  We mention here that the idea in the proof of Lemma \ref{l:m2} can be used to give an alternative proof for Theorem 6 in \cite{N16}.

We next construct a graph that is a variant of $M_k(n)$.
For the graph $M_k(n)$, we recall three subsets of vertices $X$, $Y$, and $Z$. For $k \geq 3$ and $n \geq 2k+1$, let $M'_k(n)$ be a graph obtained from $M_k(n)$ by connecting two vertices from $X$ and deleting  two edges with endpoints from $Z$.  We observe  $M'_k(n)$ is Hamiltonian.  One can easily show $\lambda(M'_k(n))<n-k-1$ and  $q(M'_k(n)) \geq 2(n-k-1)$ for $n$ large enough. Therefore,  we are not able to use Theorem \ref{Ni16} to determine whether it contains a Hamiltonian cycle. However,   it satisfies the condition in Theorem \ref{t:thm2} and we can tell it is Hamiltonian.  Similarly, we construct a graph that is a  variant of $B_k(n)$.  For $k \geq 3$ and $n \geq 2k$, recall  subsets $X, Y, Z,$ and $W$ of the vertex set of $B_k(n)$.  Starting from $B_k(n)$, we connect two vertices in $X$ and delete two edges  from $E(Y \cup W,Z)$.  Let $B'_k(n)$ be the resulting graph. It is not hard to check $\lambda(B'_k(n))< \lambda(B_k(n))$. Thus we are not able to  use Theorem \ref{lin1} to determine the Hamiltonicity of $B'_k(n)$. Since one can check $q(B'_k(n)) \geq 2n-k$ easily,   Theorem \ref{t:thm2} implies that it is Hamiltonian.  In conclusion, Theorem \ref{t:thm1} and Theorem \ref{t:thm2} provide new power for one to  determine the Hamiltonicity of graphs with large minimum degree.

\end{document}